\documentclass[12pt]{article}
\usepackage[dvips]{graphicx}
\usepackage{amsmath,amssymb}
\usepackage{cases}
\newcommand{\bm}[1]{\mbox{\boldmath $#1$}}

\newcommand{\qed}{ \hfill $\square$ \\ }
\def\Proof{\noindent{\bf Proof.}\quad}

\newtheorem{theorem}{Theorem}[section]
\newtheorem{lemma}[theorem]{Lemma}
\newtheorem{corollary}[theorem]{Corollary}

\newtheorem{example}[theorem]{Example}

\begin{document}
\title{Permutation test for dendrograms and its application to the analysis of mental lexicons}

\author{Kei Kobayashi
(The Institute of Statistical Mathematics, Tokyo) 
\and
Mitsuru Orita (Kumamoto University)}%

\maketitle
\begin{abstract}
A novel type of permutation tests for dendrogram data is studied with respect to two types of metrics for measuring the difference between dendrograms. First, the Frobenius norm is used, and we prove the consistency and efficiency of the permutation tests.  Next, the geodesic distance on a dendrogram space is used. The uniqueness of the geodesics on every dendrogram space is proved and some existing algorithms for computing geodesics are applied. Mental lexicons of English words are analyzed as an application example of the proposed permutation tests. The difference of mental lexicons 
between native and non-native English speakers is examined by analyzing sorting 
task data that used English words taken from various word classes.
\footnote[0]{This work was supported by JSPS KAKENHI Grant Numbers 24700288 and 25370634 and ISM Research Collaboration Grant 25-Kyoken-2067.}
\end{abstract}

\section{Introductory Remarks}
A dendrogram is a tree diagram usually used for representing a hierarchical clustering of a set of observed samples.
Dendrograms are used in various academic areas of data analysis including statistics, computational biology, 
psychology, and machine learning \cite{jain1999}\cite{rencher2012}\cite{webb2003}.
In the paper, we propose a novel method for the statistical hypothesis testing of the difference of dendrograms.
In order to illustrate our method effectively, sorting task data is considered 
throughout the paper, but the method can be applied to the
analysis of dendrograms in various disciplines.

\begin{example}[Experimental data: sorting English words by meaning]
\label{ex:sorting}
\rm
In our experiments, we give a task to two participant groups, native English speakers (NS) and non-native English speakers (all Japanese, therefore denoted as JP):
``Sort the given English word cards into groups of words that you think would go together according to meaning.'' 
\end{example}

The main purpose of the experiments is to analyze the sorting task results representing NS and JS differences in the English mental lexicon. 
Throughout the paper, we assume that a mental lexicon, irrespective of whether it is a mental lexicon of an individual person or an ``average'' mental lexicon of a group of people, is modeled by a dendrogram, which is a tree structure. 
Our specific goal is to analyze the NS and JP difference in dendrograms of an``average'' mental lexicon. We will omit ``average'' from the term and call it a mental lexicon for short. 
Tree models for a mental lexicon have been studied in psycholinguistic areas \cite{miller1969}\cite{rapoport1972}\cite{routh1994}. 
Of course, there are other options for a model of a mental lexicon including directed graphs 
\cite{degroot1993}\cite{kroll2005}\cite{meara2002}\cite{sanchez2004}\cite{sunderman2006}
and ``cobweb-like'' networks \cite{aitchison2003}\cite{haastrup2000}\cite{wilks2007}\cite{wolter2006}. 
Though model selection between them is an interesting  and challenging problem, we will focus on the case of a dendrogram model in this paper.

In section \ref{sec:LW-method}, we explain how to estimate a mental lexicon
from the experimental data of word sort.
We use an algorithm called the Lance-Williams method, which is
popularly used for hierarchical clustering analysis.

In section \ref{sec:test-Frobenius}, permutation test statistics under a null hypothesis of the equivalence
of two dendrograms is proposed. Furthermore, sufficient conditions for consistency and efficiency of the permutation
test are proved in section \ref{sec:theory}.

For computing the test statistics, we need to introduce a distance (metric) 
measuring the difference of two tree diagrams. 
We first use the Frobenius norm between two distance matrices
computed using the path length in each dendrogram.
Meanwhile, in section \ref{sec:test-geodesic}, we use a geodesic distance on the set of 
dendrograms (dendrogram space).
A dendrogram space is a subset of a tree space, which has been studied recently, 
in particular, for phylogenetic tree analysis \cite{billera2001}\cite{holmes2005}\cite{nye2011}\cite{miller2012}.
From the geodesic-convexity of a dendrogram space in a tree space, algorithms for 
tree spaces can also be applied to dendrogram spaces.

In section \ref{sec:experiment}, we explain the results of permutation tests 
that are applied to sorting task data probing into NS and JS differences in
mental lexicons. The differences are tested by using the two different distances 
and the results are compared.

Note that there are two existing studies with some similarity to our proposed method;
however we considered our problem independently.
First, in \cite{arnaoudova2010}, a hypothetical testing for measuring congruence of two 
phylogenetic trees via Bayesian estimation is proposed. 
They used the Markov chain Monte Carlo (MCMC) method to obtain samples from a posterior distribution of trees. 
Our method is simpler than theirs because in our setting, we can permute the 
data of two groups and compare the dendrograms more easily than for
phylogenetic trees in general. 
Even if the permutation approach can be undertaken for phylogenetic trees under 
some specific situations, it is difficult to verify the method theoretically.
Meanwhile, as we will see in section \ref{sec:theory}, some assymptotic properties
of the permutation test for dendrograms can be proved by the local linearlity of
dendrogram construction.

Second, the method proposed in section \ref{sec:test-geodesic} 
has some similarity to the algorithm proposed in \cite{chakerian2012}
in the sense that theory and algorithms for phylogenetic trees are applied
to the hierarchical clustering analysis.
Our method and theory stated in section \ref{sec:test-geodesic}
are different from their study in \cite{chakerian2012} 
in the following points: (1) we propose permutation tests while they 
mainly target the confidence interval of an estimated tree by bootstrapping, and
(2) a dendrogram space as a subset of tree space is studied and
only geodesics in the dendrogram space are used for defining the geodesic distance.

Since a geodesic in a dendrogram space is also a geodesic in a tree space embedding
as proved in section \ref{sec:theory-geodesic}, the same algorithms for computing
the geodesics can be applied to the both spaces.
However, it is important to note the difference between the definitions of 
the two spaces when we consider the motivation to introduce the geodesic distance for measuring the 
difference of trees instead of other measures.

\section{Algorithm for computing dendrograms: Lance-Williams method}
\label{sec:LW-method}

In this section, we summarize how to construct a dendrogram of
a mental lexicon from the experimental data of the sorting tasks of English words.
We use the Lance-Williams method, which is one of the most popular methods for
hierarchical clustering \cite{lance1967}.

Let $M$ be the total number of English words. 
If $n$ out of $N$ examinees classify words $W_i$ and $W_j$ ($i,j=1,\dots,M$) in the same group,
define $d_H(i,j):=1-n/N$ as a distance between $W_i$ and $W_j$. This $d_H(\cdot,\cdot)$ can be
recognized as a Hamming distance, and therefore, it satisfies the axiom of distance.
Let $\mathcal{P}_M$ be the set of partitions of the index set $\{1,\dots,M\}$
of the English words $\{W_i\}$. Therefore, each element of $\mathcal{P}_M$
corresponds to a clustering of the words.
The Lance-Williams method is popular as a method for hierarchical clustering,
but it can be recognized as a transform of a distance:
it computes another distance $d_T$ from the Hamming distance $d_H$
by the algorithm in Table \ref{table:LW-algorithm}.


\begin{table}[htb]
\caption{Lance-Williams algorithm}
\begin{center}
\fbox{
\begin{minipage}[c]{15cm}
\begin{enumerate}
\setlength{\leftskip}{1.5cm}
\item[\bf INPUT:] $d_H(i,j)$ for $i,j=1,\dots,M$.
\item[\bf STEP\;0:] Let $\mathcal{G}:=\{\{1\},\dots,\{M\}\}\in \mathcal{P}_M$ and 
$d(\{i\},\{j\}):=d_H(i,j)$ for $i,j=1,\dots,M$.
\item[\bf STEP\;1:] Select one pair $I,J\in \mathcal{G}$ attaining the minimum value of $d(I,J)$
either randomly or deterministically. Remove $I$ and $J$ from $\mathcal{G}$ and add $I\cup J$ to $\mathcal{G}$ instead. 
\item[\bf STEP\;2:] With constants $\alpha_I, \alpha_J, \beta,$ and $\gamma$ defined by Table \ref{table:LW}, set
\begin{align}
d(I\cup J,&K)=d(K,I\cup J)\nonumber\\
 &:=\alpha_I d(I,K) +\alpha_J d(J,K) + \beta d(I,J) + \gamma |d(I,K)-d(J,K)|.
\label{eq:LW} 
 \end{align}
\item[\bf STEP\;3:] 
For each $i\in I\cup J$ and $j\in K$, set $d_T(i,j):=d(I\cup J,K)$.
\item[\bf STEP\;4:]
Repeat STEP 1-3 until $\mathcal{G}$ becomes $\{\{1,\dots, M\}\}$, where $\{\{1,\dots, M\}\}$ is a class whose component is only the whole set. 
\item[\bf OUTPUT:] $d_T(i,j)$ for $i,j=1,\dots,M$.
\end{enumerate}
\end{minipage}
}
\label{table:LW-algorithm}
\end{center}
\end{table}


\begin{table}[htb]
\caption{Parameter values for the Lance-Williams method. ($n_I$ is the size of $I$.)}
\begin{center}
\begin{tabular}{|l|c|c|c|c|}
\hline
 & $\alpha_I$ & $\alpha_J$ & $\beta$ & $\gamma$ \\ \hline
Group average method & $n_I/n_{I\cup J}$ & $n_J/n_{I\cup J}$ & 0 & 0 \\ \hline
Centroid method & $n_I/n_{I\cup J}$ & $n_J/n_{I\cup J}$ & $-n_I n_J/n_{I\cup J}^2$ & 0 \\ \hline
\parbox[c][1.1cm][c]{3cm}{Ward method} & $\displaystyle \frac{n_I+n_K}{n_{I\cup J} +n_K}$ & $\displaystyle \frac{n_J+n_K}{n_{I\cup J} +n_K}$ 
 & $\displaystyle \frac{n_K}{n_{I\cup J} +n_K}$ & 0 \\ \hline
Nearest neighbor method & 1/2 & 1/2 & 0 & -1/2 \\ \hline
Furthest neighbor method & 1/2 & 1/2 & 0 & 1/2 \\ \hline
\end{tabular}
\label{table:LW}
\end{center}
\end{table}

The Lance-Williams method includes various clustering methods
obtained by setting the value of parameters $\alpha_I$, $\alpha_J$, $\beta$, and $\gamma$ 
as in Table \ref{table:LW}.
For example, the group average method defines the distance between
two clusters by the average distance of all pairs of elements (words) 
from each of the two clusters.

\begin{figure}[tb]
 \begin{center}
   \includegraphics[height=9cm,width=15cm]{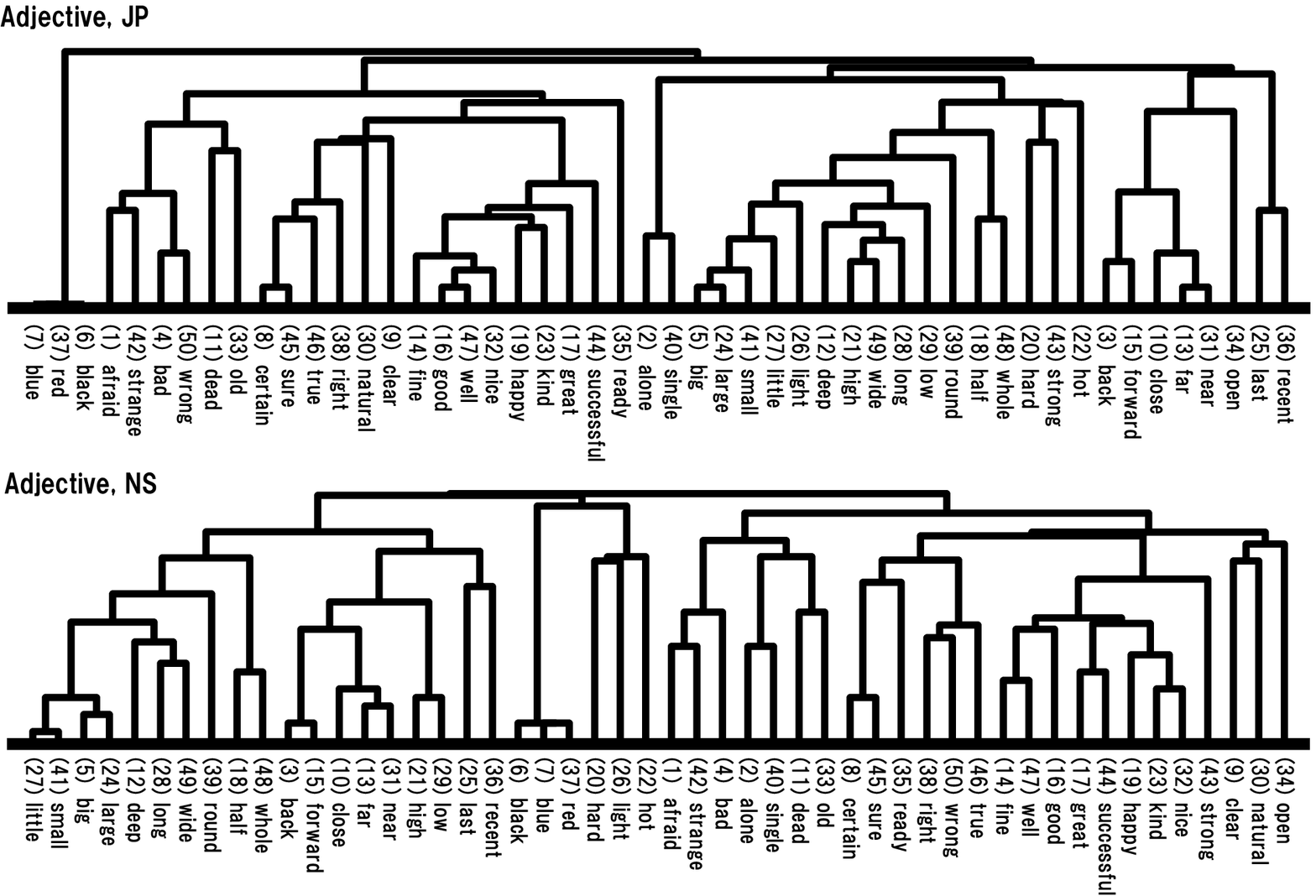}
 \end{center}
 \caption{Dendrograms of mental lexicons of English adjectives (Top: non-native speakers, Bottom: native speakers)}
 \label{fig:dendrogram_adj}
\end{figure}

For each output of the Lance-Williams method, we can construct a dendrogram
as in Figure \ref{fig:dendrogram_adj}. The dendrogram is obtained as follows:
(i) initially locate all leaf nodes $p_i$ for $i=1,\dots,M$ on the base line 
(ii) in each STEP 1 of the Lance-Williams algorithm,
connect two nodes corresponding to two clusters $I,J\in \mathcal{G}$ 
by ``$\Pi$-shaped'' line segments whose height from the base line is $d(I,J)/2$.
Relocate some of the leaf nodes if it is necessary, 
(iii) create new node $p_{I\cup J}$ corresponding to the cluster $I\cup J$
at the center of the horizontal line of the``$\Pi$-shape,''
(iv) repeat (i)-(iii) until the Lance-Williams algorithm finishes; and
(v) finally, rescale the height of the entire diagram to 1.

Note that, the order of the leaf nodes on the base line in the final dendrogram is not unique.
However, if we release all leaf nodes from the base line and contract each horizontal line of $\Pi$
to the center point, a metric tree (or a non-negatively weighted semi-labeled rooted tree) 
is obtained by retaining the length (or weight) of the edges corresponding to the vertical line-segments.
The metric tree is unique, and we identify this metric tree with the original dendrogram
throughout the paper.

For hierarchical clustering, clusters can be obtained by the sets of 
the elements (words) whose distance between each other is less than a constant $c>0$,
which is usually set heuristically. This clustering can be recognized as a separation
of the dendrogram by cutting it at a height $\tilde{c}$; a rescaling of $c$.
Meanwhile, in this paper, we consider the dendrogram itself as a representation of a mental lexicon.

Before closing this section, we remark that the Lance-Williams algorithm is not 
necessarily a projection of a distance, i.e., if we input $d_H:=d_T$
and run the same algorithm again, then the output can be different from $d_T$.
However, the group average method, the nearest neighbor method, and the furthest neighbor method
are projections by the following lemma.
\begin{lemma}
\label{lem:projection}
If $\alpha_I+\alpha_J=1$ and $\beta=0$ in the Lance-Williams algorithm,
the algorithm becomes a projection of a distance.
\end{lemma}
This is because, in (\ref{eq:LW}) of STEP 2, $d(I,K)=d(J,K)$ implies $d(I\cup J,K)=d(I,K)=d(J,K)$.
On the contrary, the centroid method and the Ward method are not projections
since $d(I,K)=d(J,K)$ implies $d(I\cup J,K)<d(I,K)$ for the controid method and
$d(I\cup J,K)>d(I,K)$ for the Ward method.
By this fact, these two methods are not preferable for our usage since 
we are estimating a dendrogram by the algorithm.

\section{Permutation test for dendrograms using the Frobenius norm}
\label{sec:test-Frobenius}

The main objective of this study is to propose a statistical hypothesis test to 
identify two dendrograms. We consider the example of sorting English words 
and denote the words by $W_1,W_2,\dots,W_M$ and the two groups of the examiniees
by GP1 and GP2. We assume, for simplicity, that the number $N$ of the examiniees in each group
is equal and an even number though this assumption is inessential and can be removed by
a slight modification.

Let $T_{\rm GP1}$ be the distance matrix whose $(i,j)$-th element is $d_T(i,j)$ computed by 
the Lance-Williams method with a specific value of $\alpha_I$, $\alpha_J$, $\beta$, and $\gamma$
for GP1. $T_{\rm GP2}$ for GP2 is defined in the same manner.
Note that $T_{\rm GP1}=T_{\rm GP2}$ if and only if the two dendrograms coincide as metric trees.
Therefore, we consider a hypothesis test:
\begin{equation}
H_0:~~T_{\rm GP1}= T_{\rm GP2},~~~~ H_1:~~T_{\rm GP1}\neq T_{\rm GP2}.
\label{Test}
\end{equation}

In order to define a test statistic, we need to set a distance to measure the difference between
$T_{\rm GP1}$ and $T_{\rm GP2}$.
Here, we use the Frobenius norm $d_F(T,T')=\|T-T'\|=(\displaystyle\sum_{i,j} |T_{ij}-T'_{ij}|^2)^{1/2}$,
the natural distance between two matrices.
Let ${\rm GP}_\sigma$ be a randomly generated group of $N$ examinees composed by
the random sampling of $N/2$ people out of each GP1 and GP2 without replacement.
Denote a group of the remaining $N$ examinees by ${\rm GP}_{\bar{\sigma}}$.
Thus, a permutation test statistic $S_{\rm p}$ is defined as
\begin{equation}
S_{\rm p} := P(\|T_{\sigma}-T_{\bar{\sigma}}\|>\|T_{\rm GP1}-T_{\rm GP2}\|),
\label{S_p}
\end{equation}
where $T_\sigma$ and $T_{\bar{\sigma}}$ are dendrograms computed using
the Lance-Williams method for groups ${\rm GP}_\sigma$ and ${\rm GP}_{\bar{\sigma}}$, respectively.

For real data analysis, we use an empirical version $\hat{S}_{\rm p}$ instead, which is
computed by (i) generating $K$ i.i.d. samples $\sigma_1,\dots,\sigma_K\sim\sigma$ and
(ii) computing the ratio of $\sigma_i$s satisfying $\|T_{\sigma_i}-T_{\bar{\sigma}_i}\|>\|T_{\rm GP1}-T_{\rm GP2}\|$.

Similar to an ordinary binomial proportion confidence interval \cite{brown2001},
a confidence interval for $\hat{S}_{\rm p}$ with the error percentile $\alpha$ is given by
$$\hat{S}_{\rm p}\pm  z_\alpha\sqrt{\hat{S}_{\rm p}(1-\hat{S}_{\rm p})/K}$$
where $z_\alpha$ is $1-\alpha/2$ percentile of a standard normal distribution.
The Wilson score interval \cite{wilson1927} yields an improved interval
\begin{equation}
\frac{2K\hat{S}_{\rm p}+z_\alpha^2 \pm z_\alpha \sqrt{4K\hat{S}_{\rm p}(1-\hat{S}_{\rm p})+z_\alpha^2}}{2(K+z_\alpha^2)},
\label{wilson-score}
\end{equation}
These intervals can be used to check if the repetition times $K$ are sufficient for 
the permutation test.

\section{Theoretical validity of the permutation test}
\label{sec:theory}

Let ${\bm X}^{(i)}\in \{0,1\}^{M\times M}$ for $i=1,\dots,N$
be a random symmetric matrix whose components ${\bm X}_{jk}^{(i)}$ 
are 0 if the $i$-th person classifies the words $w_j$ and $w_k$ in the same class
and 1 otherwise.
Then, ${\bm D}:=\frac{1}{N}\sum_{i=1}^N {\bm X}^{(i)}$ becomes a Hamming distance matrix.
From a data of sorting the words, we consider general computation methods for a dendrogram
and do not restrict to the Lance-Williams method in this section.
Such computation of a dendrogram can be recognized as a map $\tau$ between distance matrices as
$$\tau: {\bm D}\mapsto {\bm T},$$
where ${\bm T}$ is the distance matrix of the tree metric of the dendrogram.

Instead of the distance matrices themselves, we will consider random vectors
 $X^{(i)}$, $D$, and $T \in \mathbb{R}^{M(M-1)/2}$ with the components in the upper triangle
of a distance matrix ${\bm X}^{(i)}$, ${\bm D}$, and  ${\bm T}$, respectively.
A notation $\tau$ is used for a map $D\mapsto T$ besides a map $ {\bm D}\mapsto {\bm T}$.
We assume $X^{(i)}$ of the groups GP1 and GP2 are i.i.d. sampled from
distributions $P_1$ and $P_2$ with the means $\mu_1$ and $\mu_2$ and 
the nondegenerate covariances $\Sigma_1$ and $\Sigma_2$, respectively
$$
X_1^{(1)},\dots,X_1^{(N)}~\mathop{\sim}^{\rm i.i.d.}~ P_1~\mbox{and}~
X_2^{(1)},\dots,X_2^{(N)}~\mathop{\sim}^{\rm i.i.d.}~ P_2,~\mbox{independently and}
$$
$$D_1=\bar{X_1}~\mbox{and}~D_2=\bar{X_2}.$$

Next we define a distance vector $D$ of 
a group generated by a random permutation of the samples.
For simplicity, the sample number $N$ is assumed to be even.
Let $\sigma$ be a random vector uniformly taking a value in
$\{1,2\}^{N}$ such that
half of the components are $1$ and the other half are $2$, and
let $\bar{\sigma}$ be a random vector whose components $\bar{\sigma}(i)=1$ if $\sigma(i)=2$
and $\bar{\sigma}(i)=2$ if $\sigma(i)=1$.
Define $D_\sigma=\sum_{i=1}^N X_{\sigma(i)}^{(i)}/N$ and 
$D_{\bar{\sigma}}=\sum_{i=1}^N X_{\bar{\sigma}(i)}^{(i)}/N$.

Given $X_1$ and $X_2$, let $(T_1(1),T_2(1)),\dots,(T_1(K),T_2(K))$ be
independent samples distributed identically to $(T_\sigma,T_{\bar\sigma})$.
Then the permutation test statistic is
$$\hat{S}_{\rm p} = \hat{S}_{\rm p}(X_1,X_2) := \bigl|\{j=1,\dots,K \mid \|T_1(j)-T_2(j)\|>\|T_1-T_2\|\}\bigr|/K.$$
As $K$ increases, by the law of large numbers, $\hat{S}_{\rm p}$ converges to the expectation
$$S_{\rm p} = S_{\rm p}(X_1,X_2):= P_\sigma(\|T_\sigma-T_{\bar\sigma}\|>\|T_1-T_2\|).$$

Since dendrogram computation $\tau$ is generally not continuous even for the most popular algorithms
such as the group average method, we introduce a condition on $\tau$ 
in order to prove the convergence of $S_p$.
We say that $\tau$ is {\it $(D_1,D_2)$-distinguishable} if $\tau(D_1)\neq \tau(D_2)$ and
there is a constant $c>0$ such that for any $\tilde{D}_1$ and $\tilde{D}_2$ sufficiently
close to $D_1$ and $D_2$, respectively, $\|\tau(\tilde{D}_1)-\tau(\tilde{D}_2)\|>c$.
Remark that if $\tau(D_1)\neq \tau(D_2)$ and $\tau$ is continuous at $D_1$ and $D_2$, $\tau$ is $(D_1,D_2)$-distinguishable.

\begin{theorem}
\label{thm:permutation}
(1) If $\tau$ is continuous at $(\mu_1+\mu_2)/2$ and 
$(\mu_1,\mu_2)$-distinguishable, the permutation test is consistent:
$$S_{\rm p}\mathop{\rightarrow}^{\rm a.s.} 0 ~~\mbox{as}~N \rightarrow \infty.$$
(2) Furthermore, assume $\tau$ is locally linear at $(\mu_1+\mu_2)/2$ and
continuous at $\mu_1$ and $\mu_2$.
Let $\tau((\mu_1+\mu_2)/2)=A(\mu_1+\mu_2)/2$ with a matrix $A$, then for any $\epsilon>0$,
\begin{align*}
\lim_{N \rightarrow \infty}{\rm Prob}\Bigl(1-&F_{\chi^2}(\|\tau(\mu_1)-\tau(\mu_2)\|+\epsilon; A(\Sigma_1+\Sigma_2)A^\top/N) \\
&\leq S_{\rm p}
\leq 1-F_{\chi^2}(\|\tau(\mu_1)-\tau(\mu_2)\|-\epsilon; A(\Sigma_1+\Sigma_2)A^\top/N)
\Bigr) = 1
\end{align*}
as $N\rightarrow \infty$ where $F_{\chi^2}(\cdot; \Sigma)$ is the cumulative distribution function of the generalized 
$\chi^2$ distribution with a covariance matrix parameter $\Sigma$. 
\end{theorem}
\Proof
We first prove (1).
The random variables $D_\sigma$ and $D_{\bar\sigma}$ are i.i.d. to  
$N^{-1} \sum_{i=1}^{N/2} Z_{1i} + Z_{2i}$
where
$Z_{1i}~\displaystyle\mathop{\sim}^{\rm i.i.d.}~ P_1~\mbox{and}~
Z_{2i} ~\displaystyle\mathop{\sim}^{\rm i.i.d.}~ P_2.$
By the central limit theorem,
\begin{align*}
\sqrt{N}\{D_\sigma-(\mu_1+\mu_2)/2\} &\mathop{\rightarrow}^{\rm d}
N(0,(\Sigma_1+\Sigma_2)/2)~\mbox{and}~\\
\sqrt{N}\{D_{\bar\sigma}-(\mu_1+\mu_2)/2\} &\mathop{\rightarrow}^{\rm d}
N(0,(\Sigma_1+\Sigma_2)/2)~\mbox{as}~N\rightarrow \infty.
\end{align*}
Since $D_\sigma$ and $D_{\bar{\sigma}}$ are independent,
$$\sqrt{N}(D_\sigma-D_{\bar\sigma}) \mathop{\rightarrow}^{\rm d}
N(0,\Sigma_1+\Sigma_2).$$
Similarly, by the law of large numbers, $D_\sigma-D_{\bar\sigma}\displaystyle\mathop{\rightarrow}^{\rm a.s.} 0$.
Therefore, by the continuity of $\tau$ at $(\mu_1+\mu_2)/2$,
$\|\tau(D_\sigma)-\tau(D_{\bar\sigma})\|
\displaystyle\mathop{\rightarrow}^{\rm a.s.} 0$.

Meanwhile, since $D_1\displaystyle\mathop{\rightarrow}^{\rm a.s.} \mu_1$ and 
$D_2\displaystyle\mathop{\rightarrow}^{\rm a.s.} \mu_2$, $\displaystyle\mathop{\lim\inf}_{N\rightarrow \infty} \|\tau(D_1)-\tau(D_2)\|>0$ almost surely
by the $(\mu_1,\mu_2)$-distinguishable property of $\tau$.
Thus, 
$\displaystyle\mathop{\lim\sup}_{N\rightarrow \infty} \|\tau(D_\sigma)-\tau(D_{\bar\sigma})\|-\|\tau(D_1)-\tau(D_2)\|< 0$ almost surely and
(1) is proved.

Next, we prove (2). 
The probability of $T_\sigma$ and $T_{\bar\sigma}$ being in a neighborhood of
$T^*:=(\mu_1+\mu_2)/2$ converges to 1. Thus we will assume it.
By the assumption of the local linearity of $\tau$ at $T^*$, the distance matrices $T_\sigma$ and $T_{\bar{\sigma}}$
of dendrograms made by $D_\sigma$ and $D_{\bar\sigma}$ are represented by
$$T_\sigma=AD_\sigma~\mbox{and}~T_{\bar\sigma}=AD_{\bar\sigma},$$
respectively, with a matrix $A\in \mathbb{R}^{d(d-1)/2\times d(d-1)/2}$.
Therefore,
$$\sqrt{N}(T_\sigma-T_{\bar\sigma}) \mathop{\rightarrow}^{\rm d}
N(0,A(\Sigma_1+\Sigma_2)A^\top)$$
and
$$N\|T_\sigma -T_{\bar\sigma}\|_2^2 \mathop{\rightarrow}^{\rm d}
\chi^2(A(\Sigma_1+\Sigma_2)A^\top),$$
the generalized $\chi^2$ distribution with a covariance matrix parameter $A(\Sigma_1+\Sigma_2)A^\top$.

Meanwhile, $\|D_1-D_2\|=\|\mu_1-\mu_2\|+o_P(1)$ and by the assumption of the continuity of
$\tau$ at $\mu_1$ and $\mu_2$, $\|T_1-T_2\|=\|\tau(\mu_1)-\tau(\mu_2)\|+o_P(1)$.
Therefore the theorem holds.
 \qed

Now consider the property of a map $\tau$ of the Laurence-Williams algorithm
(Table \ref{table:LW-algorithm}). 
Let $\mathcal{D}$ and $\mathcal{T}$ be the set of the distance vectors $D$ and 
the set of the tree distance vectors $T$, respectively. 
Then $\mathcal{T}$ is a polyhedral subset of a polyhedral cone $\mathcal{D}\subset \mathbb{R}^{M(M-1)/2}$.

Assume that no tie occurs in the Lance-Williams algorithm
and the algorithm becomes deterministic.
Then, the topology of an output dendrogram is determined by which pair $I,J\in \mathcal{G}$
attains the minimum value of $d(I,J)$ in STEP 1. Since each $d(I,J)$ is 
defined by a combination of linear maps and absolute value operations
in (\ref{eq:LW}) of STEP 2, 
the preimage of the trees with a same topology is
defined by a set of strict inequalities and becomes the interior of a polyhedron.
Moreover, it is easy to see each preimage becomes a polyhedral cone.
Therefore $\mathcal{D}$ can be divided into a polyhedral complex (or a polyhedral fan)
denoted by $\mathcal{CD}_\tau$ 
such that every interior point of a facet of $\mathcal{CD}_\tau$ is 
mapped to a tree with a same topology. 
Meanwhile, a point on the boundary of a facet
corresponds to a tie case and the image is determined by the rule for managing ties.

\begin{corollary}
\label{cor:LW}
For the Lance-Williams method $\tau$, 
if $\tau(\mu_1)\neq \tau(\mu_2)$ and $\mu_1$, $\mu_2$, and $(\mu_1+\mu_2)/2$ are
not on the boundary of facets in $\mathcal{CD}_\tau$, 
the permutation test becomes consistent.
Furthermore if $\gamma=0$, which includes the group average method, the centroid method,
and the Ward method,  then the assertion of Theorem \ref{thm:permutation} (2) holds.
\end{corollary}
\Proof
By the algorithm of the Lance-Williams method, $\tau$ is continuous on the interior 
of a facet in $\mathcal{CD}_\tau$ and the consistency follows.
If $(\mu_1+\mu_2)/2$ is not on the boundary of a facet in 
$\mathcal{CD}_\tau$, all distance matrices sufficiently close to
$(\mu_1+\mu_2)/2$ map to trees with a same tree-topology.
This means that the recursive computations (\ref{eq:LW}) in
STEP 2 of the Lance-Williams algorithm for each of those trees are same
and become linear if $\gamma=0$.
Since the composition of those linear computations becomes linear,
$\tau$ becomes locally linear at $(\mu_1+\mu_2)/2$. \qed


The following example shows that the assumption of the regularity at $\mu_1$, $\mu_2$, and $(\mu_1+\mu_2)/2$ is essential
in Corollary \ref{cor:LW}.
\begin{example}
We use a deterministic group average method that selects a group lexicographically 
when a tie occurs.
Let the original distance matrices be
$$D_1:=
\begin{bmatrix}
0 & 2 & 3\\
* & 0 & 2\\
* & * & 0
\end{bmatrix}
~\mbox{and}~~
D_2:=
\begin{bmatrix}
0 & 3 & 2\\
* & 0 & 2\\
* & * & 0
\end{bmatrix}
$$
Then, the computed dendrograms are
$$T_1:=
\begin{bmatrix}
0 & 2 & 2.5\\
* & 0 & 2.5\\
* & * & 0
\end{bmatrix}
~\mbox{and}~~
T_2:=
\begin{bmatrix}
0 & 2.5 & 2\\
* & 0 & 2.5\\
* & * & 0
\end{bmatrix},
$$
respectively.
However, for sufficiently small $\epsilon>0$,
$$D_1^\epsilon:=
\begin{bmatrix}
0 & 2 & 3\\
* & 0 & 2-\epsilon\\
* & * & 0
\end{bmatrix}
~\mbox{and}~~
D_2^\epsilon:=
\begin{bmatrix}
0 & 3 & 2\\
* & 0 & 2-\epsilon\\
* & * & 0
\end{bmatrix},
$$
then the computed dendrograms are
$$T_1^\epsilon:=
\begin{bmatrix}
0 & 2.5 & 2.5\\
* & 0 & 2-\epsilon\\
* & * & 0
\end{bmatrix}
~\mbox{and}~~
T_2^\epsilon:=
\begin{bmatrix}
0 & 2.5 & 2.5\\
* & 0 & 2-\epsilon\\
* & * & 0
\end{bmatrix}.
$$
Therefore, $T_1^\epsilon=T_2^\epsilon$ though $T_1\neq T_2$,
and this can cause the inconsistency of the permutation test.
Here, we considered a deterministic group average method for simplicity,
but similar problems occur even if we use random selections for ties.
\end{example}


\section{Permutation test for dendrograms by the geodesic distance}
\label{sec:test-geodesic}

The permutation test proposed in section \ref{sec:test-Frobenius} uses the Frobenius norm for
measuring the difference between two dendrograms.
The middle point of two distance matrices $D_1$ and $D_2$ (i.e., matrices whose elements $D_{ij}=d(i,j)$ 
satisfy the axiom of distances) in the sense of the Frobenius norm is
$(D_1+D_2)/2$, which is also a distance matrix.
However, the middle point of two tree matrices $T_1$ and $T_2$ (i.e., matrices whose element $T_{ij}=d_T(i,j)$
is the path length between two leaves $i$ and $j$ in a non-negatively weighted rooted tree diagram)
is not necessarily a tree matrix.

In this sense, we can say that the Frobenius norm measures the difference of the distances defined by the dendrograms
rather than a difference of the dendrograms themselves.
In this section, we propose a permutation test for dendrograms with a geodesic distance 
defined on the set of dendrograms. 
Owing to the use of the geodesic distance, not only the middle point of any two dendrograms, 
but also any points on the shortest path between two dendrograms becomes a dendrogram.

\subsection{Tree space}
In this section, we define a {\it tree space} proposed in \cite{billera2001}.
In a dendrogram, the depth of each leaf is one, i.e., the path length from each leaf vertex
to the root is equal to one and the height of the dendrogram is one.
However, first, discard this condition and consider the set $\mathcal{T}_p$ of the non-negatively 
weighted rooted tree diagrams with $p$ labeled leaves. For a dendrogram of the English words,
$p$ is equal to $M$, the number of the words.

Denote the leaves and their index set as $l_1,\dots,l_p$ and $\mathcal{I}_p=\{1,\dots,p\}$, respectively.
An edge having a leaf $l_i$ as its endpoint is called as a leaf edge and is denoted by $e_{\{i\}}$.
Each edge that is not a leaf edge is called as an inner edge.
By removing an inner edge from a tree, the tree is separated into two subtrees
and only one of them includes the root node.
The removed inner edge is denoted by $e_A$, where
$A\subset \mathcal{I}_p$ is the index set of the leaves in the subtree that
does not include the root node.
Let $d_A$ be the length (or weight) of an edge $e_A$ and 
set $d_A=0$ if $e_A$ does not exist.

For example, if $p=3$, only one of $e_{\{1,2\}}$, $e_{\{1,3\}}$ or $e_{\{2,3\}}$ can exist
and only one of them can have a positive length.
In general, the following statement holds for each weighted rooted tree:
\begin{center}
(C1) for any non-trivial index sets $A, B\subset \mathcal{I}_p$,\\
$d_A>0~\mbox{and}~d_B>0 \Rightarrow A\subset B, B \subset A~\mbox{or}~A\cap B=\phi.$
\end{center}
On the other hand, if (C1) holds, there exists a unique weighted rooted tree whose edge lengths
are equal to $d_A$ for every non-trivial $A\subset \mathcal{I}_p$.

Let $v_p$ be a vector with an arbitrary order of the components
$\{d_A\}$ for non-trivial $A\subset \mathcal{I}_p$.
Since $\mathcal{I}_p$ can be non-trivially bipartitioned in $2^{p-1}-1$ different ways, 
$v_p$ becomes a $(p+2^{p-1}-1)$-dimensional vector.
Thus, the set $\mathcal{T}_p$ of the trees can be embedded into $\mathbb{E}^{p+2^{p-1}-1}$,
the $(p+2^{p-1}-1)$-dimensional Euclidean space. 
For example, if $p=3$, 
$$\mathcal{T}_3=\{(d_{\{1\}}, d_{\{2\}}, d_{\{3\}}, d_{\{1,2\}}, d_{\{1,3\}}, d_{\{2,3\}})\in \mathbb{R}_{\geq 0}^6
\mid d_{\{i\}^c}=0 \mbox{ for at least two of } i=1,2,3 \} $$
and $\mathcal{T}_3$ is embedded in $\mathbb{E}^6$.
The embedded $\mathcal{T}_p$ becomes a simplicial complex called a simplicial fan.

From the embedding, a geodesic distance on $\mathcal{T}_p$ is introduced naturally by
the shortest path length (in the sense of the Euclidean distance) among the paths through $\mathcal{T}_p$.
We call $\mathcal{T}_p$ with the geodesic distance as a {\it tree space} with $p$ leaves,
and it is denoted by $\mathcal{T}_p$.
Note that $\mathcal{T}_p$ is uniquely defined independently of the order of the components of $v_p$
up to the isometric equivalence.

\begin{figure}[tb]
 \begin{center}
   \includegraphics[height=7cm]{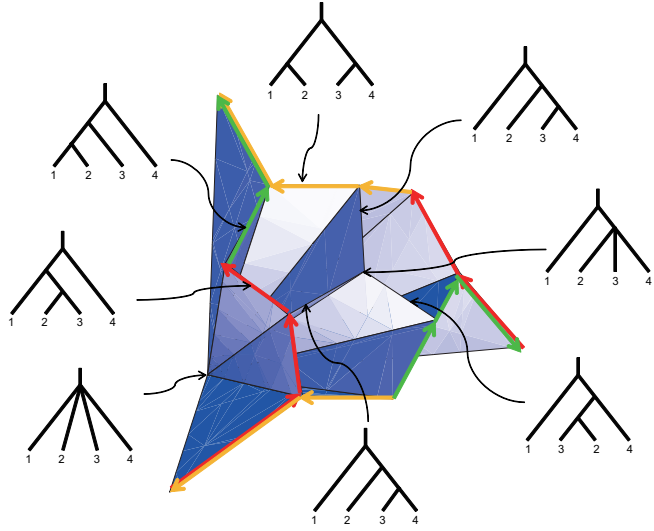}
 \end{center}
 \caption{Sketch of the subspace of a tree space $\mathcal{T}_4$ corresponding to the inner edges. This figure is inspired 
by some figures in \cite{billera2001}.}
 \label{fig:tree-space}
\end{figure}

A subspace of a tree space $\mathcal{T}_4$ corresponding to the inner edges is sketched in Figure \ref{fig:tree-space}.
The whole $\mathcal{T}_4$ becomes the direct product of the subspace and an orthant $\mathcal{R}_{\geq 0}^4$.
Since the subspace is embedded in $\mathbb{E}^7$, here we used some ``tricks'' to depict it in three dimension:
(1) axes for orthants are not orthogonal, (2) each colored cycle must be closed but is cut in the figure, and 
(3) each infinite orthant is cut to be finite. 
The topologies of trees corresponding to several points on $\mathcal{T}_4$ are also illustrated.

Tree spaces have been studied recently especially in phylogenetic tree analysis.
The theoretical study of tree spaces began from \cite{billera2001}, in which they proved 
the CAT(0) property of a tree space and the uniqueness of the geodesic between each pair of the points.
By using this property, statistical inferences and hypothetical testing on a tree space
were studied in \cite{holmes2005}.
After being applied in an innovative polynomial time algorithm in \cite{owen2011},
it has been applied to a variety of statistical inferences including the computation 
of the Fr\'{e}chet mean
\cite{bacak2012}\cite{miller2012}  and principal component analysis \cite{nye2011}.

\subsection{Dendrogram space and the permutation test by the geodesic distance}

In a dendrogram, the depth of each leaf is one. 
Let $\tilde{\mathcal{T}}_p$ be the set of the rooted trees whose depth of each leaf is one.
Therefore each tree in $\tilde{\mathcal{T}}_p$ is uniquely represented by
the set of edge lengths $d_{A:A^c}$ for non-trivial $A\subset \mathcal{I}_p$
satisfying (C1) and
\begin{center}
(C2) for any $i \in \mathcal{I}_p$,
$\displaystyle\sum_{i \in A\subsetneq \mathcal{I}_p} d_A=1.$
\end{center}

Since $\tilde{\mathcal{T}}_p$ is a subset of $\mathcal{T}_p\subset \mathbb{E}^{p+2^{p-1}-1}$, 
a geodesic distance is introduced by the shortest path length among the paths
through only in $\tilde{\mathcal{T}}_p$.
We call $\tilde{\mathcal{T}}_p$ with the geodesic distance as a {\it dendrogram space} with
$p$ leaves, and it is denoted by $\tilde{\mathcal{T}}_p$.

It is not evident that a geodesic in $\tilde{\mathcal{T}}_p$ 
is also a geodesic in $\mathcal{T}_p$.
In section \ref{sec:theory-geodesic}, we prove that this statement is true,
which is equivalent to the fact that a dendrogram space satisfies a property
called CAT(0). This is good from the computational aspect
since existing rapid algorithms for computing the geodesic 
distance in a tree space are also available for a dendrogram space.

The geodesic distance between two dendrograms $T_1,T_2\in\tilde{\mathcal{T}}_p$
is denoted by $d_{\rm geo}(T_1,T_2)$.
Then, the permutation test statistics by the geodesic distance $d_{\rm geo}$ becomes
\begin{equation}
S_{\rm geo} := P_\sigma(d_{\rm geo}(T_{\sigma},T_{\bar{\sigma}})>d_{\rm geo}(T_{\rm GP1},T_{\rm GP2})),
\label{S_geo}
\end{equation}
instead of (\ref{S_p}) by the Frobenius norm and its empirical version $\hat{S}_{\rm geo}$.

In this section, we have considered only dendrograms whose height is one.
We denote such dendrograms by DG0.
The assumption on the height is required because of
the normalization of a dendrogram computed by the Lance-Williams algorithm.
If we remove the normalization and consider a higher value on the scaling of the
original distance $d_0$ (not necessarily the Hamming distance in general), 
the assumption on the height for a dendrogram space must be removed.
If $d_0$ is a Hamming distance and takes a value in $[0,1]$
but there is no normalization after the Lance-Williams algorithm,
the height of the dendrogram is at most $1$ because of the property
of Lance-Williams method.
Therefore, we may consider two other types of dendrograms:
\begin{itemize}
\setlength{\leftskip}{1cm}
\item[(DG1)]  dendrograms whose hight can be arbitrarily large and
\item[(DG2)]  dendrograms whose hight is at most 1.
\end{itemize}

If we embed the dendrogram spaces naturally in a Euclidean space,
DG1 becomes an infinite cone ($v\in \tilde{\mathcal{T}}_p \Rightarrow c v \in \tilde{\mathcal{T}}_p$ for $c\geq 0$)
and DG2 becomes a finite cone ($v\in \tilde{\mathcal{T}}_p \Rightarrow c v \in \tilde{\mathcal{T}}_p$ for $c\in [0,1]$, but $c v\notin \tilde{\mathcal{T}}_p$ for a sufficiently large $c$);
however, DG0 is neither of them.
Therefore, it is not evident that the results on geodesics for DG0 can apply to DG1 and DG2.
Fortunately, in section \ref{sec:theory-geodesic}, we will prove the CAT(0) property of DG1 and DG2
and the algorithms to compute geodesic distances for a tree space that 
can be directly applied to DG1 and DG2.

\subsection{Theoretical verification of the permutation test using the geodesic distance}
We can prove some asymptotic properties of the permutation test statistics $S_{\rm geo}$ using 
the geodesic distance in a similar way for $S_{\rm p}$ using the Frobenius distance.
We used two different ways for representing each tree (or dendrogram):
(i) a $(p(p-1)/2)$-dimensional vector of the path length between each pair of leaves and
(ii) a $(p+2^{p-1}-1)$-dimensional vector of the edge length used for defining a tree space.
For avoiding confusion, we use an alphabet $W$ for the notation of the later vector.

If we fix a topology of a tree, the path length between two leaves is obtained by 
summing up the length of the edges on the path and $T$ becomes a linear transform 
of $W$.  On the other hand, a map from $T$ to $W$ is a bijection, and therefore,
also becomes linear. 
Denote the linear map by $W=BT$ with a matrix $B\in \mathbb{R}^{(p+2^{p-1}-1)\times (p(p-1)/2)}$. 

\begin{lemma}
\label{lem:geodesic-bound}
Let $W_i=BT_i$ for $i=1,2$. Then
$$ \|W_1-W_2\|\leq d_{\rm geo}(T_1,T_2) \leq \sqrt{2} \|W_1-W_2\|.$$
\end{lemma}
Lemma \ref{lem:geodesic-bound} will be proved after Lemma \ref{lem:orthants}.
By Lemma \ref{lem:geodesic-bound}, we can prove some asymptotic properties of
the permutation test using the geodesic distance.

\begin{theorem}
\label{thm:permutation-geo}
(1) If $\tau$ is continuous at $(\mu_1+\mu_2)/2$ and 
$(\mu_1,\mu_2)$-distinguishable, the permutation test is consistent:
$$S_{\rm geo}\mathop{\rightarrow}^{\rm a.s.} 0 ~~\mbox{as}~N \rightarrow \infty.$$
(2) Furthermore, assume $\tau$ is locally linear at $(\mu_1+\mu_2)/2$ and
continuous at $\mu_1$ and $\mu_2$.
Let $\tau((\mu_1+\mu_2)/2)=A(\mu_1+\mu_2)/2$ with a matrix $A$, then for any $\epsilon>0$,
\begin{align*}
\lim_{N\rightarrow \infty} {\rm Prob}\Bigl(1-&F_{\chi^2}\bigl(d_{\rm geo}(\tau(\mu_1),\tau(\mu_2))+\epsilon; 
B A(\Sigma_1+\Sigma_2)A^\top B^\top/N\bigr) \\
&\leq S_{\rm geo}
\leq 1-F_{\chi^2}\bigl(d_{\rm geo}(\tau(\mu_1),\tau(\mu_2))-\epsilon; B A(\Sigma_1+\Sigma_2)A^\top B^\top/N\bigr)
\Bigr) = 1
\end{align*}
\end{theorem}
\Proof
By Lemma \ref{lem:geodesic-bound},
$\displaystyle\lim_{N\rightarrow \infty} \|T_\sigma-T_{\bar{\sigma}}\|=0$ implies
$\displaystyle\lim_{N\rightarrow \infty} d_{\rm geo}(T_\sigma,T_{\bar{\sigma}})=0$.
Therefore the consistency (1) follows from Theorem \ref{thm:permutation}(1).
If we assume the local linearity of $\tau$ at $(\mu+\mu_2)/2$,
both $W_\sigma=BT_\sigma$ and $W_{\bar{\sigma}}=BT_{\bar{\sigma}}$ converge to $B(A(\mu+\mu_2)/2)$ in probability, and therefore,
the probability of $W_\sigma$ and $W_{\bar{\sigma}}$ being on a same facet of the tree space converges to one.
Thus we can assume that $d_{\rm geo}(T_\sigma,T_{\bar{\sigma}})=\|W_\sigma-W_{\bar{\sigma}}\|=\|B(T_\sigma-T_{\bar{\sigma}})\|$
and Theorem \ref{thm:permutation-geo}(2) is proved in the similar way for Theorem \ref{thm:permutation}(2).\qed

\section{CAT(0) property of the dendrogram space}
\label{sec:theory-geodesic}

In this section, we prove the CAT(0) property of dendrogram spaces.
Let $(\mathcal{X},d)$ be a geodesic metric space and denote
the geodesic (shortest path) between $a$ and $b$ in $\mathcal{X}$
by $\widetilde{ab}$. In particular for the Euclidean space, the geodesic line segment
is denoted by $\overline{ab}$.
A geodesic metric space $\mathcal{X}$ is a {\it CAT(0) space} iff
for any $a,b,c\in \mathcal{X}$ satisfies the following CAT(0) property:
``Construct a triangle in $\mathbb{E}^2$ with vertices $a', b' ,c'$, called the comparison triangle, 
such that
$\|a'-b'\|=d(a,b)$, etc. Select $p\in \widetilde{bc}$ and find the corresponding point 
$p'\in \overline{b'c'}$ such that $d(b,p) = \|b'-p'\|$.
Then for any choice of $p\in \widetilde{bc}$, $d(p,a) \leq \|p'-a'\|.$''
Intuitively speaking, each geodesic triangle in $\mathcal{X}$ is ``thinner'' than the corresponding one in a Euclidean space. 
It is known that a CAT(0) space has a nonpositive curvature locally. 

The CAT(0) space and its generalization CAT(k) space was proposed and 
studied by M. Gromov \cite{gromov1987} and a tree space was 
proved to be CAT(0) in \cite{billera2001}.
From this fact, we can prove that the geodesic between two points becomes unique
and the Fr\'echet mean is uniquely defined on a tree space.

We will prove that a dendrogram space also becomes CAT(0).  
We consider dendrograms with height 1 unless there is a further remark.

\begin{lemma}
\label{lem:orthants}
The dendrogram space $\tilde{\mathcal{T}}_m$ is a polyhedral complex such that
each facet $\tilde{\Lambda}$ is a subset of $(m-2)$-dimensional orthant $\Lambda$.
Furthermore, (i) $\tilde{\Lambda}$ is convex, (ii) the normal projection of $\tilde{\Lambda}$
to each $(m-3)$-dimensional sub-orthant of $\Lambda$ stays in $\tilde{\Lambda}$
and (iii) for each pair of adjacent facets $\tilde{\Lambda}\subset \Lambda$ and $\tilde{\Lambda}'\subset \Lambda'$ of $\tilde{\mathcal{T}}_m$, we can embed $\tilde{\Lambda}$ and a reflection of 
$\tilde{\Lambda}'$ to an $(m-1)$-dimensional Euclidean space $\mathbb{E}^{m-1}$ such that
$\Lambda$ and ${\rm Ref}(\Lambda')$ correspond to 
$\{x_1,\dots,x_{m-1}\geq 0\}$ and $\{x_1\leq 0,x_2,\dots,x_{m-1}\geq 0\}$, respectively, and
$\tilde{\Lambda}\cup {\rm Ref}(\tilde{\Lambda}')$ is convex in $\mathbb{E}^{m-1}$.
\end{lemma}
\Proof
Properties (i) and (ii) are evident from the style of the inequalities defining $\tilde{\Lambda}$.
For proving the convexity in (iii), let $p$ and $q$ be a pair of points in 
$\tilde{\Lambda}\cup {\rm Ref}(\tilde{\Lambda}')\subset \mathbb{E}^{m-1}$. If both $p$ and $q$ are in
either $\tilde{\Lambda}$ or ${\rm Ref}(\tilde{\Lambda}')$, their convex combination
is evidently in $\tilde{\Lambda}$ or ${\rm Ref}(\tilde{\Lambda}')$, respectively, and therefore, 
we assume 
$p\in \tilde{\Lambda}\cap {\rm Ref}(\tilde{\Lambda}')^c$ and $q\in {\rm Ref}(\tilde{\Lambda}') \cap \tilde{\Lambda}^c$.
Let $p_0$ and $q_0$ be the normal projections of $p$ and $q$, respectively, to the 
$(m-2)$-dimensional orthant $\Lambda\cap {\rm Ref}(\Lambda')$. Then, by (ii) of the lemma,
$p_0, q_0 \in \tilde{\Lambda}\cap {\rm Ref}(\tilde{\Lambda}')$.
Next, set $r$ to be the point internally dividing the line segment $\overline{p_0q_0}$ in the ratio $\|\overline{pp_0}\|$ and $\|\overline{qq_0}\|$. Then, $r$ is also in $\tilde{\Lambda}\cap {\rm Ref}(\tilde{\Lambda}')$
by (i) of the lemma.
Now, the line segments $\overline{pr}$ and $\overline{rq}$ are in 
$\tilde{\Lambda}$ and ${\rm Ref}(\tilde{\Lambda}')$, respectively, and the points $p$, $r$, and $q$ are on a single line.
This means that the convex combination of $p$ and $q$ is in $\tilde{\Lambda}\cup {\rm Ref}(\tilde{\Lambda}')$
and the convexity of $\tilde{\Lambda}\cup {\rm Ref}(\tilde{\Lambda}')$ follows. \qed

Before proving CAT(0) property of the dendrogram space, 
we will prove Lemma \ref{lem:geodesic-bound}.\\

(Proof of Lemma \ref{lem:geodesic-bound}) 
The first inequality is evident since $d_{\rm geo}(T_1,T_2)$ is a path length between 
$T_1,T_2\in\mathcal{T}_p\subset \mathbb{E}^{p+2^{p-1}-1}$ while $\|W_1-W_2\|$ is the shortest path
length in $\mathbb{E}^{p+2^{p-1}-1}$. 
For showing the second inequality, let $F_1$ (or $F_2$) be a facet
including $T_1$ (or $T_2$, respectively) in the dendrogram space.
If $F_1$ and $F_2$ are not adjacent, $W_1\perp W_2$ and
$$d_{\rm geo}(T_1,T_2) \leq \|W_1\|+\|W_2\|\leq \sqrt{2}(\|W_1\|^2+\|W_2\|^2)^{1/2} =\sqrt{2}\|W_1-W_2\|.$$
Thus we assume that $F_1$ and $F_2$ are adjacent.
Consider the orthogonal projections of $T_1$ and $T_2$ to a face 
$F_1\cap F_2$ and denote them by $P_1$ and $P_2$, respectively. 
By setting $p:=T_1$, $q:=T_2$, $p_0:=P_1$, and $q_0:=P_2$, we can apply the notation used
in the proof of Lemma \ref{lem:orthants}(iii).
Since $\|\overline{pq}\|^2=(\|\overline{pp_0}\|+\|\overline{qq_0}\|)^2+\|\overline{p_0q_0}\|^2$,
\begin{align*}
d_{\rm geo}(T_1,T_2)^2&=(\|W_1-P_1\|+\|W_2-P_2\|)^2+\|P_1-P_2\|^2\\
&\leq  2(\|W_1-P_1\|^2+\|W_2-P_2\|^2) + 2\|P_1-P_2\|^2\\
&=2\|W_1-W_2\|^2.
\end{align*}
Here we used that $W_1-P_1$, $W_2-P_2$, and $P_1-P_2$ are orthogonal to each other. \qed

\begin{theorem}
\label{thm:CAT0}
The dendrogram space $\tilde{\mathcal{T}}$ is CAT(0).
\end{theorem}

\Proof
Let $\mathcal{T}$ be the tree space that is the cone generated by $\tilde{\mathcal{T}}$.
Since every tree space $\mathcal{T}$ is CAT(0),
it is sufficient to prove that for every $x,y\in \tilde{\mathcal{T}}$
the geodesic on $\mathcal{T}$ between $x$ and $y$ 
is also a geodesic on $\tilde{\mathcal{T}}$.

Assume there is a geodesic on $\tilde{\mathcal{T}}$
that is not a geodesic on $\mathcal{T}$. Then we can select a pair of points
$p,q\in \partial\tilde{\mathcal{T}}$ such that the whole geodesic $\overline{pq}$ of $p,q$ 
on $\mathcal{T}$ connecting $p,q$ excepting the two ends lies outside $\tilde{\mathcal{T}}$. 
Here, $\partial\tilde{\mathcal{T}}$ means the boundary of $\tilde{\mathcal{T}}$
as a subset of $\mathcal{T}$.

Now, the cone ${\rm cone}(\widetilde{pq})$ generated by the geodesic 
$\widetilde{pq}$ is a (two-dimensional) flat surface
and isometrically embeddable in $\mathbb{E}^2$ such that the origin of the original space is
embedded as the origin of $\mathbb{E}^2$. Note that the length of the geodesic
$|\widetilde{pq}|$ is shorter than $\|\overline{pO}\|$ and $\|\overline{qO}\|$, since
otherwise, the geodesic $\overline{pq}$ is a sequence of $\overline{pO}$ and $\overline{Oq}$
and the geodesic on $\mathcal{T}$ and the geodesic on $\tilde{\mathcal{T}}$ must coincide.
Therefore, $\widetilde{pq}$, $\overline{pO}$, and $\overline{qO}$ are isometrically
embedded in $\mathbb{E}^2$
as line segments $\overline{p'q'}$, $\overline{p'O'}$, and $\overline{q'O'}$
for $p',q', O'\in\mathbb{E}^2$, respectively.

The corresponding embedding in $\mathbb{E}^2$
of ${\rm cone}(\widetilde{pq}) \cap \partial\tilde{\mathcal{T}}$
becomes a sequence of line segments connecting $p'$ and $q'$.
Denote the sequence $p'=p'_0,p'_1,\dots,p'_m=q'$ and 
consider $\zeta_i:= \langle (\overrightarrow{O'p'}+\overrightarrow{O'q'})/2, \overrightarrow{O'p'_i}\rangle
/\|(\overrightarrow{O'p'}+\overrightarrow{O'q'})/2\| \|\overrightarrow{O'p'_i}\|$.
Then, $\zeta_0<0$ and $\zeta_{m-1}>0$,
and therefore, there exists $1\leq i\leq m-2$ satisfying $\zeta_i<\zeta_{i+1}$.
If the corresponding $p_i$ in the original space $\tilde{\mathcal{T}}$ is an interior point
of a facet orthant, we can select points $a'\in \overline{p'_{i-1}p'_i}$ and 
$b'\in \overline{p'_{i}p'_{i+1}}$, which are so close to $p'_i$ that
the corresponding points $p,a,b$ are in the same orthant.
However, because $\zeta_i<\zeta_{i+1}$, the shortest path from $a$ to $b$
must lie outside $\tilde{\mathcal{T}}$ and this contradicts the convexity of the
facet of $\tilde{\mathcal{T}}$ (Lemma \ref{lem:orthants} (ii)).

If $p_i$ is on an intersection of two facets, we can select points $a'\in \overline{p'_{i-1}p'_i}$ and 
$b'\in \overline{p'_{i}p'_{i+1}}$, which are so close to $p'_i$ that
the corresponding points $a$ and $b$ are in each of the two adjacent facets.
Even if we embed these two facets to a Euclidean space as in Lemma \ref{lem:orthants} (iii),
by $\zeta_i<\zeta_{i+1}$ the shortest path from $a$ to $b$
must lie outside $\tilde{\mathcal{T}}$. This contradicts the convexity of 
the union of the two facets in Lemma \ref{lem:orthants} (iii).\qed

By the proof of Theorem \ref{thm:CAT0}, 
algorithms to compute geodesic distances for a tree space $\mathcal{T}$
can be used for a corresponding dendrogram space $\tilde{\mathcal{T}}$
since every geodesic on $\tilde{\mathcal{T}}$ is also a geodesic on $\mathcal{T}$.

In section \ref{sec:test-geodesic}, we introduced DG1, the set of the dendrograms
whose hight can be arbitrarily large, and DG2, that whose hight is at most 1.
The CAT(0) property of DG1 can be derived from that of DG2.
This is because every geodesic triangle for checking the CAT(0) property
on a dendrogram space $\mathcal{T}_{\rm DG1}$ of  $DG1$  
can be included in a dendrogram space DG2($c$) of the dendrograms whose hight is at most $c$
for a sufficiently large $c$. Since DG2($c$) is just a rescaling of DG2,
CAT(0) property of DG2 implies that of DG2($c$) and therefore DG1.

By condition (C2) defining a dendrogram space of DG0,
$\sum_A d_A\leq 1$, where the summation is over non-trivial 
$A\subset\mathcal{I}_p$ such that $i\in A$ and $|A|\geq 2$.
A dendrogram space of DG2 is defined by modifying the inequalities 
as $\sum_A d_A+\tilde{w}\leq 1$ with a non-negative variable $\tilde{w}$
corresponding to the length of a shortest leaf edge.
Then, each facet is a subset of a $(p-1)$-dimensional orthant, 
which is one dimension higher than a facet of the dendrogram space DG0.
We can check easily that each facet is convex and 
the normal projection of each facet to each $(p-2)$-dimensional sub-orthant stays
in the facet. By the same arguments for Lemma \ref{lem:orthants}(iii) and
Theorem \ref{thm:CAT0}, the CAT(0) property of DG2 is proved.

\section{Experimental results on sorting task data}
\label{sec:experiment}

We computed the permutation test statistics $\hat{S}_{\rm p}$ for 
experimental data of sorting tasks of Example \ref{ex:sorting}.
Five different sets of English words, all of which are arbitrarily selected from 
the JACET List of 8000 Basic Words \cite{JACET}, are used: 
(i) the most frequent 500 words (0.5K) of mixed word classes; 
(ii) the first, 1000 high frequent 50 verbs (${\rm 1K}_\alpha$); and (iii) the first, 1000 high frequent 50 verbs ((iv) adjectives or (v) nouns) (${\rm 1K}_\beta$), 
where the last three sets are carefully chosen for stronger semantic links between words.
The number of participants in each group, NS or JP, is 28 for the word set (i) and 30 for the other four word sets. 
See \cite{orita2011} for the details of the experiments. 

We used the group averaging method to make a dendrogram.
The number of generated random permutations is set $K=5000$ and
the Frobenius norm is used.
A Matlab program ``linkage.m'' in the Statistics Toolbox is used for the group average method.
We list the values of $\hat{S}_{\rm p}$ in Table \ref{table:P-values}.
A conservative correction of $\hat{S}_{\rm p}$ by the Wilson score interval (\ref{wilson-score})
with a significance level 1\% is denoted by $\hat{S}_{\rm p}^*$ and also listed.

\begin{table}[t h]
 \caption{Value of the permutation test statistics $\hat{S}_{\rm p}$ for data (i)-(v).}
 \begin{center}
 \begin{tabular}{|l || r|r|r|r|r|}
 \hline
 Data & (i) & (ii) & (iii) & (iv) & (v) \\
 \hline
 Word set &  mixed &  verbs  & verbs  & nouns  & adjectives  \\
 & (0.5K) & (${\rm 1K}_\alpha$) & (${\rm 1K}_\beta$) & (${\rm 1K}_\beta$) & (${\rm 1K}_\beta$) \\
 \hline
 \#Participants & 30 & 28 & 30 & 30 & 30 \\
 \hline
 $\hat{S}_{\rm p}$ (\%) & 1.22 & 49.60 & 0.42 & 2.18 & 0.04 \\
 \hline
 $\hat{S}_{\rm p}^*$  (\%) & 1.53 & 51.02 & 0.63 & 2.58 & 0.17 \\
 \hline
 \end{tabular}
 \label{table:P-values} 
\end{center} 
\end{table}

If we consider $\hat{S}_{\rm p}$ as the p-value of the permutation test,
the null hypothesis $T_{\rm NS} = T_{\rm JP}$ is
rejected with a significant level 5\% for data (i), (iii), (iv), and (v),
whereas the null hypothesis is rejected with a significant level 1\% only for data (iii) and (v).
Since the corrected $\hat{S}_{\rm p}^*$ implies the same result,
we can conclude that the number of permutations, 5000, is sufficiently large.
Besides for testing the hypothesis, we can use the statistics $\hat{S}_{\rm p}$ as
a relative measure for comparing the difference between $T_{\rm NS}$ and $T_{\rm JP}$.

Next, we compare the results obtained through
the permutation test statistics $\hat{S}_{\rm p}$ with the Frobenius norm
and $\hat{S}_{\rm geo}$ with the geodesic distance.
The word sets (iii) verbs (${\rm 1K}_\beta$), (iv) nouns (${\rm 1K}_\beta$), and (v) adjectives (${\rm 1K}_\beta$), stated above are used.
In addition to the NS and JP sorting results, in analysis, we included the results produced by a group of 30 novice Japanese learners of English (NV) who have significantly lower English proficiency than JP. 
The number of permutations is 5000.
We used the geodesic distances computed by Geodesic Treepath Problem (GTP) algorithm 0.1
by Megan Owen and J. Scott Provan \cite{owen2011}.
The values of $\hat{S}_{\rm p}$ and $\hat{S}_{\rm geo}$ are listed in Table \ref{table:Frobenius-geodesic}.

\begin{figure}[bt]
 \begin{center}
   \includegraphics[height=7cm]{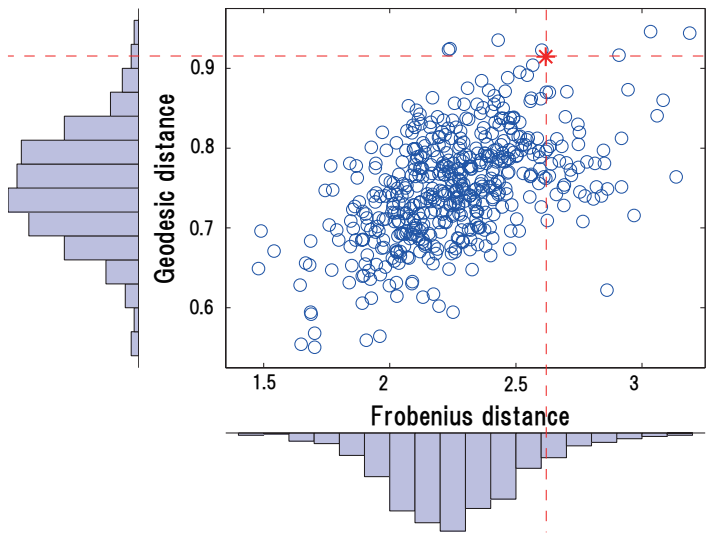}
 \end{center}
 \caption{Frobenius norm and the geodesic distance: the blue ``o''s are 
$(\|T_{\sigma}-T_{\bar{\sigma}}\|, d_{\rm geo}(T_{\sigma},T_{\bar{\sigma}}))$
for each random permutation $\sigma$ and the red ``$*$'' is 
$(\|T_{\rm GP1}-T_{\rm GP2}\|, d_{\rm geo}(T_{\rm GP1},T_{\rm GP2}))$.}
 \label{fig:euclid_geodesic}
\end{figure}

\begin{table}[t h]
 \caption{Value of the permutation test statistics $\hat{S}_{\rm p}$ and  $\hat{S}_{\rm geo}$ for data (iii)-(v).} 
\begin{center}
 \begin{tabular}{|l | l | l || r|r|}
 \hline
 GP1 & GP2 & Word set & $\hat{S}_{\rm p}$  & $\hat{S}_{\rm geo}$ \\
 \hline \hline
 NS & JP & (iii) verbs & 0.44 & 0.58\\
 \hline 
 NS & JP & (iv) nouns & 2.06 & 2.12\\
 \hline 
 NS & JP & (v) adjectives & 0.04 & 0.12\\
 \hline \hline
 JP & NV & (iii) verbs & \bm{7.66} & \bm{1.00} \\
 \hline 
 JP & NV & (iv) nouns & \bm{0.38} & \bm{1.76}\\
 \hline 
 JP & NV & (v) adjectives & 0.06 & 0.08\\ 
 \hline \hline
 NV & NS & (iii) verbs & 0.02 & 0.00\\
 \hline 
 NV & NS & (iv) nouns & 0.00 & 0.00\\
 \hline 
 NV & NS & (v) adjectives & 0.00 & 0.00\\
 \hline
 \end{tabular} 
\label{table:Frobenius-geodesic}
 \end{center} 
\end{table}
From the table, the null hypothesis $T_{\rm JP}= T_{\rm NV}$ for the data set (iii)
is rejected with a significant level 5\% by $\hat{S}_{\rm p}$, but it is not rejected by $\hat{S}_{\rm geo}$.
Figure \ref{fig:euclid_geodesic} depicts the different results yielded by the two distances.
On the contrary, the null hypothesis $T_{\rm JP}= T_{\rm NV}$ for the data set (iv)
is not rejected with a significant level 1\% by $\hat{S}_{\rm p}$, but it is rejected by $\hat{S}_{\rm geo}$.

\section{Concluding remarks and discussions}
In the paper, we proposed the permutation tests for dendrograms
with two different distances on the dendrograms: the Frobenius norm 
and the geodesic distance.
We proved some asymptotic properties of the permutation test statistics
by the Frobenius norm.
For the geodesic distance, we proved that a dendrogram space, the set of 
dendrograms with a geodesic metric naturally inherited from an embedding Euclidean space,
has CAT(0) property. 
Therefore the algorithms to compute the geodesics on a tree space, which have been developed 
and used in phylogenetic analysis, can be directly applied to a tree space.

Then, which distance should we use? This is a natural question
and the answer depends on the context of the hypothesis testing.
First, the Frobenius norm is measuring the difference between distance matrices
computed by the path lengths between the leaves in each dendrogram.
Thus, a dendrogram can be recognized as a method of approximation of the distances between
the leaves, whereas the geodesic distance measures the difference of each edge length, and
therefore, it focuses on the tree structure itself.
For example, a natural method of defining the average tree of
trees $T_1$ and $T_2$ by the Frobenius norm is $\tau((T_1+T_2)/2)$ by recognizing $T_1$ and $T_2$
as distance matrices. This depends on the tree construction $\tau$, whereas the average by
the geodesic distance, which is the middle point of the geodesic
between $T_1$ and $T_2$, is independent of $\tau$.

Another difference is their computational costs: computation of the Frobenius norm 
requires $O(p^2)$, whereas the fastest algorithm for the geodesic distance requires $O(p^4)$ \cite{owen2011}.
Therefore, if the number of samples becomes much larger, using the Frobenius norm may
be only the feasible method.


The Lance-Williams method for computing dendrograms includes various methods
such as the group average method, the centroid method, the Ward method, the nearest neighbor method,
and the furthest neighbor method.
It is worth remarking that the group average method is the most preferable in our context
since it is the only one which is a projection (Lemma \ref{lem:projection}) and also has the local linearity
for proving the asymptotic efficiency of the permutation test (Corollary \ref{cor:LW}).

For verifying the permutation testing by the asymptotic theory, 
we need the assumption that the sample size is sufficiently large.
One of the methods to check if the sample size is sufficiently large or not is 
to make a confidence interval of the permutation statistics by 
employing the bootstrapping method \cite{efron1979}.



In \cite{miller2012} and \cite{bacak2012}, the Fr\'echet mean, computed by minimizing 
the sum of the squared geodesic lengths, is studied.
As we proved in the paper, a dendrogram space has CAT(0) property and
it implies that the Fr\'echet mean of dendrograms in a dendrogram space
also becomes a dendrogram.
Therefore, instead of computing a dendrogram by the Lance-Williams method,
the Fr\'echet mean of the mental lexicons of each examinee can be a candidate
of an ``average'' mental lexicon of a group of examinees.
For this use of the Fr\'echet mean, we need to manage the computational cost since
it usually requires much more cost than the Lance-Williams algorithm.




\thispagestyle{empty}

\end{document}